\documentclass[a4paper,12pt]{amsart}
\usepackage{amssymb}
\usepackage{ifthen}
\usepackage{graphicx}
\nonstopmode \numberwithin{equation}{section}
\newtheorem{thm}[equation]{Theorem}
\newtheorem{cor}[equation]{Corollary}
\newtheorem{lem}[equation]{Lemma}
\newtheorem{prop}[equation]{Proposition}

\newtheorem{conj}{Conjecture}

\theoremstyle{definition}

\newtheorem{prob}[equation]{Problem}

\newenvironment{rem}{%
\medskip
\noindent \textsl{{\sl Remark. }}}{
}
\newenvironment{rems}{%
\bigskip
\noindent \textsl{{\sl Remarks. }}}{\bigskip}

\newcounter {own}
\def\theown {\thesection       .\arabic{own}}

\newenvironment{pf}[1][]{%
 \vskip 3mm
 \noindent
 \ifthenelse{\equal{#1}{}}%
  {{\slshape Proof. }}%
  {{\slshape #1.} }%
 }%
{\qed
\medskip
}

\newcounter{alphabet}
\newcounter{tmp}

\newcommand{\ID}{{\mathbb D}}

\newcommand{\IC}{{\mathbb C}}

\def\be{\begin{equation}}
\def\ee{\end{equation}}

\newcommand{\bee}{\begin{enumerate}}
\newcommand{\eee}{\end{enumerate}}

\newcommand{\blem}{\begin{lem}}
\newcommand{\elem}{\end{lem}}
\newcommand{\bthm}{\begin{thm}}
\newcommand{\ethm}{\end{thm}}
\newcommand{\bcor}{\begin{cor}}
\newcommand{\ecor}{\end{cor}}
\newcommand{\beg}{\begin{examp}}
\newcommand{\eeg}{\end{examp}}
\newcommand{\begs}{\begin{examples}}
\newcommand{\eegs}{\end{examples}}
\newcommand{\bdefe}{\begin{defin}}
\newcommand{\edefe}{\end{defin}}
\newcommand{\bprob}{\begin{prob}}
\newcommand{\eprob}{\end{prob}}
\newcommand{\bei}{\begin{itemize}}
\newcommand{\eei}{\end{itemize}}

\newcommand{\bcon}{\begin{conj}}
\newcommand{\econ}{\end{conj}}
\newcommand{\bcons}{\begin{conjs}}
\newcommand{\econs}{\end{conjs}}
\newcommand{\bprop}{\begin{prop}}
\newcommand{\eprop}{\end{prop}}
\newcommand{\br}{\begin{rem}}
\newcommand{\er}{\end{rem}}
\newcommand{\brs}{\begin{rems}}
\newcommand{\ers}{\end{rems}}
\newcommand{\bo}{\begin{obser}}
\newcommand{\eo}{\end{obser}}
\newcommand{\bos}{\begin{obsers}}
\newcommand{\eos}{\end{obsers}}
\newcommand{\bpf}{\begin{pf}}
\newcommand{\epf}{\end{pf}}
\newcommand{\ba}{\begin{array}}
\newcommand{\ea}{\end{array}}
\newcommand{\beq}{\begin{eqnarray}}
\newcommand{\beqq}{\begin{eqnarray*}}
\newcommand{\eeq}{\end{eqnarray}}
\newcommand{\eeqq}{\end{eqnarray*}}

\newcommand{\ds}{\displaystyle}

\begin{document}
\bibliographystyle{amsplain}
\title{On some problems of James Miller}
\author{B. Bhowmik}
\address{B. Bhowmik, Department of Mathematics,
Indian Institute of Technology Madras, Chennai-600 036, India.}
\email{ditya@iitm.ac.in; bappaditya.bhowmik@gmail.com}
\author{S. Ponnusamy}
\address{S. Ponnusamy, Department of Mathematics,
Indian Institute of Technology Madras, Chennai-600 036, India.}
\email{samy@iitm.ac.in}
\author{K.-J. Wirths}
\address{K.-J. Wirths, Institut f\"ur Analysis, TU Braunschweig,
38106 Braunschweig, Germany.}
\email{kjwirths@tu-bs.de}

\subjclass[2000]{30C45}
\keywords{Starlike, meromorphic, and Schwarz' functions,  Taylor coefficient}
\date{ 
Feb. 5, 08; File: Mil1.tex}
\begin{abstract}
We consider the class of meromorphic univalent functions having a simple pole at
$p\in(0,1)$ and that map the unit disc onto the exterior of a domain which is starlike with
respect to a point $w_0 \neq 0,\, \infty$. We denote this class of functions by
$\Sigma^*(p,w_0)$. In this paper, we  find the exact region of variability for the
second Taylor coefficient
for functions in $\Sigma^*(p,w_0)$. In view of this result we rectify some results of
James Miller.
\end{abstract}

\thanks{}

\maketitle
\pagestyle{myheadings}
\markboth{B. Bhowmik, S. Ponnusamy and K.-J. Wirths
}{Domain of variability}

\section{Introduction}
Let  $\ID := \{ z: \,  |z| < 1 \}$ be the unit disc in the complex plane $\IC$. Let
$\Sigma^*$ denote the class of functions
$$g(z)=\frac{1}{z}+d_0+d_1z+d_2z^2+\cdots
$$
which are univalent and analytic in $\ID$ except for the simple pole
at $z=0$ and map $\ID$ onto a domain whose
complement is starlike with respect to the origin.
Functions in this class is referred to as the meromorphic
starlike functions in $\ID$.
This class has been studied by Clunie \cite{clunie-59} and later an extended
version by Pommerenke \cite{pomm-63},  and many others. Another related
class of our interest is the class $S(p)$ of univalent
meromorphic  functions $f$ in $\ID$ with a simple pole at
$z=p$, $p\in (0,1)$, and with the normalization
$f(z)=z+\sum_{n=2}^{\infty}a_n(f)z^n$ for $|z|<p$.
If $f\in S(p)$ maps $\ID$ onto a domain whose complement
with respect to $\overline{\IC}$ is convex, then we call $f$
a concave function with pole at $p$ and the class of these functions
is denoted by $Co(p)$. In a recent paper,
Avkhadiev and  Wirths \cite{AW} established the region
of variability for $a_n(f)$, $n\geq 2$, $f\in Co(p)$
and as a consequence two conjectures of Livingston \cite{Living-94} in 1994 and
Avkhadiev, Pommerenke and  Wirths \cite{APW} were settled.

In this paper, we consider the class $\Sigma^*(p,w_0)$ of meromorphically
starlike functions $f$ such that $\overline\IC \setminus f(\ID)$ is a starlike set
with respect to a finite point $w_0\neq 0$ and have the standard
normalization $f(0)=0=f'(0)-1$. We now recall the following analytic
characterization for functions in $\Sigma^*(p,w_0)$.

\bigskip
\noindent
{\bf Theorem A.} {\em
$f\in\Sigma^*(p,w_0)$ if and only if there is a probability measure $\mu (\zeta)$ on
$\partial\ID=\{\zeta: |\zeta|=1\}$ so that
$$
f(z)=w_0 + \frac{pw_0}{(z-p)(1-zp)}\exp\left(\int_{\partial\ID} 2\log (1-\zeta z)d\mu (
\zeta) \right )
$$
where
$$
w_0= -\frac{1}{p+1/p-2\int_{|\zeta|=1}\zeta d\mu(\zeta)}.
$$
}

The necessary part of Theorem~A has been proved by Miller \cite{Mil-80}
while the sufficiency part has been established by Yuh Lin \cite[Theorem 1]{Lin-88}.
In \cite{Mil-72, Mil-80}, Miller discussed a numbers of properties of the
class $\Sigma^*(p,w_0)$. See also \cite{BP1, Lin-88,Lin-Owa-90} for some other
basic results  such as bounds for  $|f(z)-w_0|$

We may state an equivalent formulation of Theorem~A (see also
\cite{Lin-Owa-90}). A function $f$ is said to be in $\Sigma^*(p,w_0)$ if and
only if there exists an analytic function $P(z)$ in $\ID$ with  $P(0)=1$ and
\be\label{p5eq0}
{\rm Re}\, P(z)> 0, ~ z\in\ID,
\ee
where
\be\label{p5eq1a}
P(z)= \frac{-zf'(z)}{f(z)-w_0}-\frac{p}{z-p}+\frac{pz}{1-pz}.
\ee
We may write $P(z)$ in the following power series form
$$P(z)=1+b_1z+b_2z^2+\cdots.
$$
Also, each $f\in \Sigma^*(p,w_0)$ has  the Taylor expansion
\be\label{p5eq1}
f(z)=z+\sum_{n=2}^{\infty}a_n(f)z^n,\quad |z|<p.
\ee
To recall the next result, we need to introduce a notation.
Let $\mathcal{P}(b_1)$ denote the class of analytic functions $P(z)$ satisfying
$P(0)=1$, $P'(0)=b_1$ and ${\rm Re}\, P(z)> 0$ in $\ID.$

In 1972, Miller \cite{Mil-72} obtained estimations for the second Taylor
coefficient $a_2(f)$. Indeed, he showed that

\bigskip
\noindent
{\bf Theorem B.} {\em
If $f(z)\in \Sigma^*(p,w_0)$, then the second coefficient is given by
$$
a_2(f)= \frac{1}{2}w_0\left(b_2-p^2-\frac{1}{p^2}-\frac{1}{{w_0}^2}\right)
$$
where $b_2$ is the second coefficient of a function in  $\mathcal{P}(b_1)$, i.e.
the region of variability for $a_2(f)$ is contained in the disc
\be\label{p5eq1aa}
\left|a_2(f)+\frac{1}{2}w_0\left(p^2+\frac{1}{p^2}+\frac{1}{{w_0}^2}\right)\right|\leq |w_0|.
\ee
Further there is a $p_0,~ 0.39<p_0<0.61$, such that if $p<p_0$, then
${\rm Re}\, a_2(f) >0$.
}
\medskip

In 1980, Miller \cite[Equation (9)]{Mil-80} also proved a sharp estimate
regarding the second Taylor coefficient. In fact, he showed that
\be\label{p5eq4}
\left|a_2(f)-\frac{1+p^2}{p}-w_0\right|\leq |w_0|, \quad f \in \Sigma^*(p,w_0).
\ee

The aim of this paper is to find the region of variability for the second coefficient
$a_2(f)$ of functions in $\Sigma^*(p,w_0)$ for any fixed pair $(p, w_0)$.
Also we find the exact region of variability for $a_2(f)$ for fixed $p$,
and as a consequence of this we show that ${\rm Re}\, a_2(f) >0$ for all
values of $p\in (0,1)$ which Miller did not
seem to expect as we see in the last part of Theorem ~B.

\section{Region of Variability of Second Taylor Coefficient for functions in
$\Sigma^*(p,w_0)$}\label{p5sec1}

\bthm\label{p5thm1}
Let $f\in \Sigma^*(p,w_0)$ having the expansion {\rm (\ref{p5eq1})}. Then for
a fixed pair $(p, w_0)$, the exact region of variability of the second Taylor coefficient $a_2(f)$ is
the disc determined by the inequality
\beq\label{p5eq2aa}\nonumber
\left|a_2(f)-\left(p+\frac{1}{p}+w_0\right)+
\frac{1}{4}w_0\left(p+\frac{1}{p}+\frac{1}{{w_0}}\right)^2\right|&&\\
&\leq& |w_0| \left(1-\frac{1}{4}\left|p+\frac{1}{p}+\frac{1}{w_0}\right|^2\right).
\eeq
\ethm\bpf
The proof uses the  representation formula (\ref{p5eq0}), i.e. $f\in \Sigma^*(p,w_0)$
if and only if ${\rm Re}\,P(z)>0$ in $\ID$ with $P(0)=1$, where $P$ is given by
(\ref{p5eq1a}).
Since it is convenient to work with the class of Schwarz functions, we can write
each such $P$ as
\be\label{p5eq2}
P(z)=\frac{1+\omega(z)}{1-\omega(z)}, \quad z\in \ID,
\ee
where $\omega \colon \ID\to \ID$ is holomorphic  with $\omega (0)=0$ so that
$\omega(z)$ has the form
\be\label{p5eq1aaa}
\omega(z)= c_1z+c_2z^2+\cdots .
\ee
Using (\ref{p5eq1a}) and the power series representations
of $P(z)$ and $f(z)$, it is easy to compute
\beq\label{p5eq0a}
\left\{
\ba{lll}
b_1&=& \ds p+\frac{1}{p}+ \frac{1}{w_0},~ \mbox{and} \\
b_2&=&\ds  p^2+ \frac{1}{p^2}+ \frac{1}{{w_0}^2}+\frac{2a_2(f)}{w_0}.
\ea
\right.
\eeq
Now eliminating $w_0$ from (\ref{p5eq0a}), we get
\be\label{p5eq2a}
b_2= p^2+\frac{1}{p^2}+\left [b_1-\left(p+\frac{1}{p}\right)\right]^2
+2a_2(f)\left [b_1-\left(p +\frac{1}{p}\right)\right ].
\ee
Using the power series representations of $P(z)$ and $\omega (z)$, it follows
by comparing the coefficients of $z$ and $z^2$ on both sides that
$$b_1=2c_1 ~ \mbox{ and } ~b_2=2(c_1^2+c_2).
$$
Inserting the above two relations in (\ref{p5eq2a}), we get
$$2(c_1^2+c_2)=p^2+\frac{1}{p^2}
+\left [2c_1-\left(p+\frac{1}{p}\right)\right]^2+2a_2(f)\left [2c_1-\left(p
+\frac{1}{p}\right)\right].
$$

Now solving the above equation for $a_2(f)$, we get
\be\label{p5eq3a}
a_2(f)=\frac{1}{p}+ p\,\left(\frac{c_1^2-c_2+p^2-2c_1p}{1+p^2-2c_1p}\right).
\ee
Now, since $w_0$ and $p$ are fixed, we have  $c_1$ fixed. Hence using the well
known estimate $|c_2|\leq 1- |c_1|^2$ for unimodular bounded function $\omega (z)$,
the last equation results the following estimate
$$
\left|a_2(f)-\frac{1}{p}- p\,\left(\frac{c_1^2+p^2-2c_1p}{1+p^2-2c_1p}\right)\right|
\leq \frac{p(1-|c_1|^2)}{|1+p^2-2c_1p|}.
$$
Now, as $b_1=2c_1$, substituting $c_1= \frac{1}{2}(p+1/p+1/w_0)$
in the above equation we get the required estimate as given in (\ref{p5eq2aa}).
A point on the boundary of the disc described by (\ref{p5eq2aa}) is attained
for the unique extremal functions
given by (\ref{p5eq1a}) and (\ref{p5eq2}), where
$$\omega(z)=\frac{z(c_1+cz)}{1+\overline{c_1}cz},\quad  |c|=1.
$$
The points in the interior of the disc described in (\ref{p5eq2aa})
are attained for the same functions, but with $|c|<1$.
\epf

\br
Comparison of Theorem ~B and Theorem \ref{p5thm2} below, shows that the
exact region of variability of $a_2(f)$ found by Miller is for the case $c_1=0$ only.
A little computation reveals that both variability regions are the same
for $c_1=0$, i.e.,
$$\left|a_2(f)-\frac{1+p^2+p^4}{p(1+p^2)}\right|\leq \frac{p}{1+p^2}.
$$

This also shows that (\ref{p5eq4}) gives the precise region of variability
only for the case $c_1=0$. In all other cases, the boundaries of the
discs described by (\ref{p5eq1aa}) and (\ref{p5eq4})
have only one point in common with the disc
described by (\ref{p5eq2aa}) because, in the both cases, on the boundaries of the discs
described by (\ref{p5eq1aa}) and (\ref{p5eq4}), we need
$|b_2|=2$. Now, as $b_2= 2(c_2+{c_1}^2)$, this means that $|c_2+{c_1}^2|=1$. According
to the coefficients bounds for unimodular bounded function, this is only possible
for a unique $c_2$ if $c_1 \neq 0$.
\er

In the following theorem, we describe the  exact region of variability of
the second Taylor coefficient of $f\in \Sigma^*(p,w_0)$, where only $p$ is fixed.

\bthm\label{p5thm2}
Let $f\in \Sigma^*(p,w_0)$ having the expansion {\rm (\ref{p5eq1})} and let $p$ be fixed. Then the exact set
of variability of the second Taylor coefficient $a_2(f)$ is given by
\be\label{p5eq3}
|a_2(f)-1/p|\leq p.
\ee
\ethm
\bpf We may rewrite (\ref{p5eq3a}) as
\be\label{p5eq3aa}
a_2(f)=\frac{1}{p}+ p\,M,
\ee
where
$$
M=\frac{c_1^2-c_2+p^2-2c_1p}{1+p^2-2c_1p}.
$$

We wish to prove that $|M|\leq 1$. Since $\omega '(0)=c_1$, we have $|c_1|\leq 1$.

Now we fix $c_1\in \overline{\ID}$. Then $c_1^2-c_2$ varies in the closed disc
$$
\Delta(c_1):= \{z: |z-c_1^2|\leq 1-|c_1|^2\}.
$$
The map
$$
T(\zeta)=\frac{\zeta+p^2-2c_1p}{1+p^2-2c_1p}
$$
maps the disc $\Delta(c_1)$ onto the disc with center
$$
\frac{c_1^2+p^2-2c_1p}{1+p^2-2c_1p}
$$
and radius
$$
\frac{ 1-|c_1|^2}{|1+p^2-2c_1p|}.
$$
Therefore, in order to prove $|M|\leq 1$, it suffices to show that
$$\left|\frac{c_1^2+p^2-2c_1p}{1+p^2-2c_1p}\right|+\frac{ 1-|c_1|^2}{|1+p^2-2c_1p|}\leq 1.
$$
This is equivalent to
$$|c_1-p|^2+1-|c_1|^2={\rm Re}\,(1+p^2-2c_1p)\leq |1+p^2-2c_1p|.
$$
We see that equality is attained in the above inequality if and only if $c_1$
is real. Now for real $c_1$, we have
$$
T(\Delta(c_1))=\overline{\ID}~ \mbox{ if and only if}~\, c_1=p\, \mbox{ or}~
w_0=\frac{-p}{1-p^2}.
$$
Hence the extremal functions for the inequality (\ref{p5eq3}) are
given by (\ref{p5eq1a}) with $P(z)$ as in  (\ref{p5eq2}) with
$$\omega(z)=\frac{z(p+cz)}{1+pcz},\, |c|=1,
$$
and the points in the interior of the disc described by (\ref{p5eq3})
are attained for the same functions, but with $|c|<1$.
We observe that for real $c_1$ we can obtain $M=1$
only for $c_2 = c_1^2-1$. This results in other starlike
centers, but the extremal function is always the same, since $a_2(f) =
p+1/p$  is attained in the class $S(p)$  only for $f(z)=z/((1-zp)(1-z/p))$,
see for instance \cite{Jenkin-62}.
\epf

\br
This result ensures us that ${\rm Re\,} a_2(f) >0$ for all $p\in (0,1)$. In the article
{\rm \cite[Theorem 1]{Mil-72}}, Miller hoped for a possibility that for $p>.61$,
the real part of $a_2(f)$ may be negative.
But in view of our theorem we conclude that his hope was in vain.
\er

\br
In \cite{Mil-80}, Miller has obtained an estimate for the real part of the
third coefficient $a_3(f)$ for all $p$. However, in geometric function
theory, the classical question of finding
the exact region of variability for $a_n(f)$,  $n\geq 3$,
$f\in \Sigma^*(p,w_0)$, remains an open problem.
\er


\begin{thebibliography}{99}
\bibitem{APW} {\sc F. G. Avkhadiev, Ch. Pommerenke, and K.-J. Wirths}: On the coefficients
of concave univalent functions, \textit{Math. Nachr.} {\bf 271}(2004), 3--9.
\bibitem{AW} {\sc F. G. Avkhadiev and K.-J. Wirths}: A proof of the Livingston conjecture,
\textit{Forum math.} {\bf 19}(2007),  149--158.
\bibitem{BP1}{\sc B. Bhowmik and S. Ponnusamy}: Coefficient inequalities for
concave and meromorphically starlike univalent functions, \textit{
Ann. Polon. Math., To appear}.
\bibitem{clunie-59} {\sc J. Clunie}:
On meromorphic schlicht functions,
\textit{J. London Math. Soc.} {\bf 34}(1959), 215--216.

\bibitem{Jenkin-62}{\sc J. A Jenkins}:
On a conjecture of Goodman concerning meromorphic univalent functions,
\textit{Michigan Math. J.} {\bf 9}(1962), 25--27.

\bibitem{Lin-88} {\sc Chang Yuh Lin}:
On the representation formulas for the functions in the class
$\Sigma^*(p,w_0)$, \textit{Proc. Amer. Math. Soc.} {\bf 103}(1988), 517--520.

\bibitem{Living-94} {\sc A. E. Livingston:} Convex meromorphic mappings,
\textit{Ann. Polon Math.} {\bf 59}(3)(1994), 275--291.

\bibitem{Mil-72} {\sc J. Miller}: Starlike meromorphic functions,
\textit{Proc. Amer. Math. Soc.} {\bf 31}(1972), 446--452.
\bibitem{Mil-80} {\sc J. Miller}:
Convex and starlike meromorphic functions,
\textit{Proc. Amer. Math. Soc.} {\bf 80}(1980), 607--613.
\bibitem{pomm-63} {\sc Ch. Pommerenke:} On meromorphic starlike functions,
\textit{Pacific J. Math.} {\bf 13}(1963), 221--235.
\bibitem{Lin-Owa-90} {\sc Zhang Yulin and S. Owa}:
Some remarks on a class of meromorphic starlike functions,
\textit{Indian J. pure appl. Math.} {\bf 21}(9)(1990), 833--840.
\end{thebibliography}
\end{document}